\documentstyle[]{article}
\long\def\comment#1{}

\def\calB{{\cal B}}

\def\calH{{\cal H}}
\def\calL{{\cal L}}

\def\calS{{\cal S}}
\def\calT{{\cal T}}

\def\Txyz{T}
\def\Phixyz{\Phi}
\def\Psixyz{\Psi}
\def\tensor{\otimes}
\def\End{\mathop{\rm End}\nolimits}
\def\Complex{{\bf C}}

\def\Integer{{\bf Z}}

\def\sn{\mathop{\rm sn}\nolimits}
\def\cn{\mathop{\rm cn}\nolimits}
\def\dn{\mathop{\rm dn}\nolimits}

\def\tr{\mathop{\hbox{\rm tr}}\nolimits}

\renewcommand{\theequation}{\thesection.\arabic{equation}}

\newcommand{\appref}[1]{Appendix~\ref{#1}}
\begin{document}
\begin{titlepage}

\title{Algebraic Bethe Ansatz for XYZ Gaudin model}

\author{E. K. Sklyanin\thanks{EKS acknowledges the support of the
Department of Mathematical Sciences, the University of Tokyo,
where the most part of the work was done.}\\
St. Petersburg Branch \\ of the Steklov Mathematical Institute,\\
Fontanka 27, St. Petersburg 191011, RUSSIA \\
\\
\\
T. Takebe\thanks{On leave of absence from Department of Mathematical
Sciences, the University of Tokyo. The work of TT is supported by
Japan Society for the Promotion of Science, Postdoctral Fellowship for
Research abroad.}\\ 
Department of Mathematics,\\ 
University of California, \\ 
Berkeley, CA 94720, U.S.A.\\}

\date{25 January 1996}

\maketitle

\centerline{\sl
Dedicated to Professor~Hikosaburo~Komatsu on his 60th birthday}

\abstract{
The eigenvectors of the Hamiltionians of the XYZ Gaudin model are
constructed by means of the algebraic Bethe Ansatz. The construction
is based on the quasi-classical limit of the corresponding results for
the inhomogeneous higher spin eight vertex model.
}

\end{titlepage}
\section{Introduction}

The XYZ Gaudin model was introduced by M.~Gaudin \cite{gau:73,gau:76,gau:83}
as a quasiclassical limit of XYZ spin-1/2 chain. Gaudin noticed also that
the former model can be generalized to any values of constituing spins.
Whereas the spectrum and eigenfunctions of the XXX and XXZ variants 
of Gaudin model can easily be found via Bethe ansatz, the general case
encounters the same problems as in case of the original XYZ model.
For the spin-1/2 XYZ model a proper generalization of Bethe ansatz was 
found in \cite{baxter:72}, see also \cite{takh-fad:79}.
However, when passing to the Gaudin model,
an additional problem arises: whereas it is
easy to perform quasiclassical limit in the Bethe equations 
determining 
the spectrum of XYZ model, the expressions for the eigenvectors appear
to behave singularly. Solution to this problem is the main subject of
the present paper.

We obtain the arbitrary spin XYZ Gaudin model as a 
quasiclassical limit of the inhomogeneous higher spin 
generalization of XYZ model introduced in \cite{take:92,take:93,take:95}.
To construct the eigenvectors of XYZ model we use the algebraic version
of Bethe ansatz proposed in \cite{takh-fad:79} for spin-1/2 and 
generalized to higher spins in \cite{take:92,take:93,take:95}.
The paper is concluded with describing a way to circumvent the divergencies
in the quasiclassical limit and to get finite expressions for Bethe 
eigenvectors of Gaudin model.

\section{Inhomogeneous XYZ model}
\setcounter{equation}{0}

First let us briefly recall the inhomogeneous XYZ spin chain
model. The model is parametrized by an anisotropy parameter $\eta$
besides the elliptic modulus $\tau$, shifts of the spectral parameter
$z_n$, and a spin at each site $\ell_n$. Later we will take the limit
$\eta \to 0$ and recover the generating function of the integrals of
motion of the XYZ Gaudin model as the leading term of the transfer
matrix of this model as in \cite{skl:87}.

The inhomogeneous XYZ model is a quantum spin chain model defined on the
Hilbert space $\calH = \bigotimes_{n=1}^N V_n$, where $V_n$ is the spin
$\ell_n$ representation space $\Theta_{00}^{4\ell_n+}$ of the Sklyanin
algebra \cite{skl:82}, \cite{skl:83}. See \appref{skl-alg} for notations.

The transfer matrix of the model, $\hat t(u)$, is a linear operator on
$\calH$ defined as the trace of the monodromy matrix, $\Txyz(u)$, in the
context of the quantum inverse scattering method \cite{takh-fad:79}:
\begin{eqnarray}
    \Txyz(u)
    &=& \pmatrix{ A_N(u) & B_N(u) \cr
                  C_N(u) & D_N(u) }
    := L_N(u-z_N) \dots L_2(u-z_2) L_1(u-z_1),
\label{def:monodromy}
\\
    \hat t(u)  &=& \tr_{\Complex^2} (\Txyz(u)) = A_N(u) + D_N(u),
\label{def:transfer}
\end{eqnarray}
where the $L$ operators, $L_n(u)$, 
whose elements act non-trivially only on the $n$-th
component, $V_n$, are defined by (cf.~(\ref{def:Lxyz}))
\begin{eqnarray}
    L_n(u) 
    &=& \pmatrix{ \alpha_n(u) & \beta_n(u)\cr
                  \gamma_n(u) & \delta_n(u)}
    := \sum_{a=0}^3 W_a^L(u) \rho^{\ell_n}_n (S^a) \tensor \sigma^a,\\
    \rho^{\ell_n}_n(S^a) &=& 1 \tensor \dots \tensor 1
                               \tensor \rho^{\ell_n}(S^a) \tensor
                             1 \tensor \dots \tensor 1.
\label{def:L-n}
\end{eqnarray}

Because of the fundamental commutation relation,
\begin{equation}
    R_{12}(u_1 - u_2) \Txyz_{01}(u_1) \Txyz_{02}(u_2)
    =
    \Txyz_{02}(u_2) \Txyz_{01}(u_1) R_{12}(u_1 - u_2),
\label{fund-comm-rel}
\end{equation}
the transfer matrices commute with each other:
\begin{equation}
    [\hat t(u_1), \hat t(u_2)] = 0.
\label{comm-rel-of-t}
\end{equation}

The quantum determinant (\cite{kul-skl}) is defined by
\begin{equation}
    \Delta(u) := 
    \tr_{12} P_{12}^- \Txyz_{01}(u - \eta) \Txyz_{02}(u + \eta),
\label{def:q-det}
\end{equation}
where $P_{12}^-$ is a projector onto the space of anti-symmetric
tensors. This operator commute with all elements of the Sklyanin algebra
and acts on the Hilbert space $\calH$ as a scalar multiplication:
\begin{eqnarray}
    \Delta(u)|_{\calH} &=& \Delta_-(u-\eta) \, \Delta_+(u+\eta),
\label{q-det-val}
\\
    \Delta_{\pm}(u)    &=& \prod_{n=1}^N 
                           \frac{\theta_{11}(u-z_n \pm 2\ell_n\eta)}
                                {\theta_{11}(u-z_n)}.
\label{def:delta-pm}
\end{eqnarray}

\section{XYZ Gaudin model}
\setcounter{equation}{0}

We define the XYZ Gaudin model as a quasi-classical limit of the
inhomogeneous XYZ model defined above. Let us examine the asymptotic
behaviour of the operators in the previous section when $\eta$
tends to 0. The $L$ operator, the monodromy matrix, the transfer
matrix, the quantum determinant and the $R$ matrix are expanded as
\begin{eqnarray}
    L_n(u)    &=& 1 + 2 \eta \calL_n(u) + O(\eta^2),
\label{exp:L-n}
\\
    \Txyz(u)  &=& 1 + 2 \eta \calT(u) + \eta^2 \calT^{(2)}(u) + O(\eta^2),
\label{exp:Txyz}
\\
    \hat t(u) &=& 1 + \eta^2 \tr\calT^{(2)}(u) + O(\eta^3),
\label{exp:hat-t}
\\
    \Delta(u) &=& 1 + \eta^2 (\tr\calT^{(2)}(u) + 4 \tr\calT(u)^2)
                    + O(\eta^3),
\label{exp:Delta}
\\
    R(u)      &=& 1 - 2 \eta r(u) + O(\eta^2).
\label{exp:R}
\end{eqnarray}
Explicitly 
$$
    \calL_n (u) =
    \sum_{a = 1}^3 w_a(u) \rho^{\ell_n}_n \calS^a \tensor \sigma^a,
\qquad
    \calT(u) = \sum_{n=1}^N \calL_n (u - z_n).
$$
Here $\rho^\ell$ denotes the spin $\ell$ representations of the Lie
algebra $sl(2)$, and 
\begin{eqnarray*}
    w_1(u) &=& \frac{\sigma_{10}(u)}{\sigma(u)} =\frac{\cn}{\sn}(u)
    = \frac{\theta'_{11}}{\theta_{10}}
      \frac{\theta_{10}(u)}{\theta_{11}(u)},
    \\
    w_2(u) &=& \frac{\sigma_{00}(u)}{\sigma(u)} =\frac{\dn}{\sn}(u)
    = \frac{\theta'_{11}}{\theta_{00}}
      \frac{\theta_{00}(u)}{\theta_{11}(u)},
    \\
    w_3(u) &=& \frac{\sigma_{01}(u)}{\sigma(u)} =\frac{1}{\sn(u)}
    = \frac{\theta'_{11}}{\theta_{01}}
      \frac{\theta_{01}(u)}{\theta_{11}(u)},
\end{eqnarray*}
where $\theta_{ab}=\theta_{ab}(0)$, 
$\theta'_{11}= d/du (\theta_{11}(u))|_{u=0}$, $\sigma$ is
Weierstra{\ss}' sigma function and 
$\calS^a$ are generators of the Lie algebra $sl(2)$:
$$
    [\calS^a, \calS^b] = i \calS^c.
$$
Here $(a,b,c)$ denotes a cyclic permutation of $(1,2,3)$. The
commutation relations of the $\calL$ and the monodromy operator are
\begin{eqnarray}
    [\calL_{01}(u), \calL_{02}(v)]
    &=&
    [r(u-v), \calL_{01}(u) + \calL_{02}(v)],
\label{comm-rel:L}
\\
    {[\calT_{01}(u), \calT_{02}(v)]}
    &=&
    [r(u-v), \calT_{01}(u) + \calT_{02}(v)].
\label{comm-rel:T}
\end{eqnarray}
The {\em classical $r$ matrix} $r$ is defined by (\ref{exp:R}), or
explicitly by
$$
    r(u) = -\frac{1}{2} \sum_{a=1}^3 w_a(u) \sigma^a\tensor \sigma^a.
$$
The commutation relation (\ref{comm-rel:L}) is nothing but the
quasi-classical limit of the fundamental commutation relation
(\ref{fund-comm-rel}).

We define the {\em XYZ Gaudin model} by specifying its generating
function of integrals of motion as $\hat\tau(u) = \tr \calT(u)^2$.
It is obvious from (\ref{exp:hat-t}) and (\ref{exp:Delta}) that 
$\hat\tau(u)$ can be expressed as follows (cf.~\cite{skl:87}):
\begin{equation}
    \Delta(u) - \hat t(u) - 1 = 4 \eta^2 \hat\tau(u) + O(\eta^3).
\label{def:hat-tau}
\end{equation}
Therefore we can expect that eigenvectors of $\hat\tau(u)$ can be
constructed as a leading term of the $\eta$ expansion of eigenvectors of
$\hat t(u)$ which are found in \cite{take:92}, \cite{take:93},
\cite{take:95}. Essentially this is true, but we must be careful in taking
the limit as we will see in the next section.

Operator $\hat\tau(u)$ is explicitly written down as follows:
\begin{equation}
    \hat\tau(u)
    = \sum_{n=1}^N \wp(u-z_n) \ell_n(\ell_n+1)
      + \sum_{n=1}^N H_n \zeta(u-z_n) + H_0,
\label{tau}
\end{equation}
where $\zeta$ and $\wp$ are Weierstra{\ss}' zeta and $\wp$ functions and
\begin{eqnarray}
    H_n &=& 
    2 \sum_{m \neq n} \sum_{a=1}^3 w_a(z_n - z_m) \calS_n^a \calS_m^a, 
    \qquad \sum_{n=1}^N H_n=0,
\\
    H_0 &=& - \sum_{n=1}^N \sum_{a=1}^3
         \left(
         e_a (\calS_n^a)^2 +
         \sum_{m \neq n} w_a(z_n - z_m)
         \int_{\omega_{\bar a}/2}^{z_n - z_m + \omega_{\bar a}/2}
         \! \! \wp(u)\, du \, \calS_n^a \calS_m^a
         \right) \nonumber \\
        &=& \sum_{n=1}^N \sum_{a=1}^3 \Bigl(
            - e_a (\calS_n^a)^2 + \nonumber \\
        & & +
              \sum_{m \neq n} w_a(z_n - z_m)
              \left(
              \zeta\Bigl(z_n - z_m +\frac{\omega_{\bar a}}{2}\Bigr)
            - \zeta\Bigl(\frac{\omega_{\bar a}}{2}\Bigr)
              \right) \calS_n^a \calS_m^a
            \Bigr)
\label{def:H-n}
\end{eqnarray}
are integrals of motion. Here $\bar a$ is $1,3,2$ for $a=1,2,3$
respectively, $\omega_1 = 1$, $\omega_2 = \tau$, 
$\omega_3 = 1 + \tau$, $e_a = \wp(\omega_{\bar a}/2)$. We omit the
symbol $\rho^\ell$ for the sake of simplicity.

\section{Algebraic Bethe Ansatz for  XYZ spin chain}
\setcounter{equation}{0}

The algebraic Bethe Ansatz in the context of the quantum inverse
scattering method is applied to the inhomogeneous XYZ spin chain
model in the following way. (cf.~\cite{takh-fad:79}, \cite{take:92},
\cite{take:95}), Hereafter we assume that the total spin
$\ell_{\rm total} = \ell_1 + \cdots + \ell_N$ is equal to an integer
$M$.

First we introduce the matrix of the gauge transformation:
\begin{equation}
    M_\lambda(u) :=
      \pmatrix{
      - \theta_{01}\left(\frac{\lambda-u}{2};\frac{\tau}{2} \right) &
      - \theta_{01}\left(\frac{\lambda+u}{2};\frac{\tau}{2} \right)\cr
        \theta_{00}\left(\frac{\lambda-u}{2};\frac{\tau}{2} \right) &
        \theta_{00}\left(\frac{\lambda+u}{2};\frac{\tau}{2} \right)
      }
      \pmatrix{
        1 & 0 \cr
        0 & \theta_{11}(\lambda;\tau)^{-1}
      }.
\label{def:M}
\end{equation}
A twisted monodromy matrix is defined by means of $M_{\lambda}$:
\begin{eqnarray}
\lefteqn{
    \Txyz_{\lambda,\lambda'}(u;v) =
    \pmatrix{A_{\lambda,\lambda'}(u;v) &  B_{\lambda,\lambda'}(u;v) \cr
             C_{\lambda,\lambda'}(u;v) &  D_{\lambda,\lambda'}(u;v)}
    }
\nonumber
\\
    &:=& M_{\lambda}(u+v)^{-1} \Txyz(u) M_{\lambda'}(u+v),
\end{eqnarray}
where $v$ is a parameter. There is a vector $\omega_{\lambda}(v)$ called
a local pseudovacuum in the spin $\ell$ representation space 
$V_\ell = \Theta^{4\ell+}_{00}$:
\begin{eqnarray*}
    \omega_{\lambda}(v) = \omega_{\lambda}(v;y)
    := &\prod_{j=1}^{2\ell}&
    \theta_{10}\left(y+\frac{\lambda}2-\frac{v}2+(2j-3\ell-1)\eta\right)
    \times
    \\
    &\times&
    \theta_{10}\left(y-\frac{\lambda}2+\frac{v}2-(2j-3\ell-1)\eta\right),
\end{eqnarray*}
and a global pseudovacuum in $\calH$:
\begin{eqnarray}
\lefteqn{
    \Omega_N^\lambda(v):=
    \omega_{\lambda+4(\ell_1 + \cdots + \ell_N)\eta}(v + z_N) \tensor
    \dots
    }\nonumber
\\
    &\dots& \tensor
    \omega_{\lambda+4(\ell_1 + \cdots + \ell_n)\eta}(v + z_n) \tensor
    \dots \tensor
    \omega_{\lambda+4\ell_1\eta}(v + z_1).
\label{def:pseudo-vac}
\end{eqnarray}
As shown in \cite{takh-fad:79} and \cite{take:92}, we can construct
eigenvectors of the transfer matrix $\hat t(u)$ of the form
\begin{eqnarray}
    \Psixyz_\nu(w_1, \dots, w_M) &=& 
    \sum_{a\in\Integer} 
    e^{2\pi i \nu a \eta} 
    \Phixyz_{\lambda+ 2a\eta}(w_1, \dots, w_M),
\label{def:Psixyz}
\\
    \Phixyz_\lambda(w_1, \dots, w_M; v) &:=&
    B_{\lambda+2\eta,\lambda-2\eta}(w_1)
    B_{\lambda+4\eta,\lambda+4\eta}(w_2) \dots \nonumber \\
    &\dots&
    B_{\lambda+2M\eta,\lambda-2M\eta}(w_M)
    \Omega_N^{\lambda-2M\eta}(v),
\label{def:Phixyz}
\end{eqnarray}
where an integer $\nu$ and the parameters $w_1, \dots, w_M$ satisfy a
system of equations (Bethe equations) for $j=1, \dots, M$:
\begin{equation}
    \frac{\Delta_+(w_j)}{\Delta_-(w_j)}
    =
    e^{-4\pi i \nu \eta}
    \prod_{k=1, k\neq j}^M
    \frac{\theta_{11}(w_j - w_k + 2 \eta)}
         {\theta_{11}(w_j - w_k - 2 \eta)}.
\label{BA-eq:xyz}
\end{equation}
The eigenvalue $t(u)$ of the transfer matrix $\hat t(u)$:
\begin{equation}
    \hat t(u) \Psixyz_\nu(w_1, \dots, w_M) 
       = t(u) \Psixyz_\nu(w_1, \dots, w_M),
\end{equation}
satisfies
\begin{equation}
    t(u) q(u) = \Delta_+(u) q(u-2\eta) + \Delta_-(u) q(u+2\eta),
\label{funct-eq:xyz}
\end{equation}
where
\begin{equation}
    q(u) = e^{-\pi i \nu u} \prod_{m=1}^M \theta_{11}(u - w_n).
\label{def:q}
\end{equation}

Note that the Bethe vector $\Psixyz_\nu(w_1,\dots,w_M)$ is expressed as a
Fourier series, and the convergence of such series is not a priori
known. In fact when $\eta = r'/r$ is a rational number, this Fourier
series diverges and we must replace it by a sum with $r$ terms. Hence we
cannot take the limit $\eta \to 0$ of the above results naively to obtain
eigenvectors of the XYZ Gaudin model.
  
\section{Algebraic Bethe Ansatz for XYZ Gaudin model}
\setcounter{equation}{0}

In order to obtain an eigenvector of the XYZ Gaudin model, we bypass
the divergent series (\ref{def:Psixyz}) in the following way. At an
intermediate stage of derivation of the Bethe equations, we have the
formula:
\begin{eqnarray}
    \lefteqn{
    \hat t(u) \Phixyz_\lambda(\vec w;v) 
    =
    \Lambda(u; \vec w; \eta) \Phixyz_{\lambda-2\eta}(\vec w;v) +
    \Lambda(u; \vec w;-\eta) \Phixyz_{\lambda+2\eta}(\vec w;v) +
    }\nonumber
    \\
    &+& 
    \sum_{j=1}^M 
        \Lambda_j^\lambda(u;\vec w; \eta) 
        \Phixyz_{\lambda-2\eta}(\vec w_j(u);v) +
    \sum_{j=1}^M 
        \Lambda_j^\lambda(u;\vec w;-\eta) 
        \Phixyz_{\lambda+2\eta}(\vec w_j(u);v),
\label{t-Phi}
\end{eqnarray}
where we abbreviated $(w_1, \dots, w_M)$ to $\vec w$ and 
$(w_1, \dots, w_{j-1}, u, w_{j+1}, \dots, w_M)$ to $\vec w_j(u)$, and
$\Lambda$ and $\Lambda_j^\lambda$ are defined by:
\begin{eqnarray}
    \Lambda(u;\vec w; \eta) &=& 
    \Delta_+(u) 
    \prod_{m=1}^M \frac{\theta_{11}(u-w_m-2\eta)}{\theta_{11}(u-w_m)}
\label{def:Lambda}
\\
    \Lambda_j^\lambda(u;\vec w; \eta) &=&
    \theta_{11}(-2\eta)
    \frac
    {\theta_{11}((u - w_j)-(\lambda - 2\eta))}
    {\theta_{11}(u - w_j) \theta_{11}(\lambda - 2\eta)}
    \times
    \\
    &\times&
    \Delta_+(w_j) 
    \prod_{k=1, k\neq j}^M
    \frac{\theta_{11}(w_j - w_k - 2 \eta)}
         {\theta_{11}(w_j - w_k)}.
\label{def:Lambda-j}
\end{eqnarray}
Since all of these formulas are holomorphic around $\eta = 0$, we can
expand them as Taylor series in $\eta$. For example, as
$M_\lambda(u)$ does not depend on $\eta$ and the monodromy matrix
$\Txyz(u)$ has the expansion (\ref{exp:Txyz}), the twisted monodromy
matrix $\Txyz_{\lambda+2m\eta,\lambda-2m\eta}(u;v)$ has the expansion,
\begin{equation}
    \Txyz_{\lambda+2m\eta,\lambda-2m\eta}(u;v)
    =
    \pmatrix{1 & 0 \cr 0 & 1} + O(\eta).
\end{equation}
Therefore (1,2)-element of the twisted monodromy matrix has an
expansion as follows:
\begin{equation}
    B_{\lambda+2m\eta,\lambda-2m\eta}(u;v)
    =
    \eta \calB_{\lambda, m}(u,v) + O(\eta^2).
\label{exp:B}
\end{equation}
The pseudovacuum vector defined by (\ref{def:pseudo-vac}) is expanded as
\begin{equation}
    \Omega_N^{\lambda-2M\eta}(v)=
           \Omega_N^{\lambda,M,(0)}(v) 
    + \eta \Omega_N^{\lambda,M,(1)}(v) +  \dots.
\label{exp:pseudo-vac}
\end{equation}
Combining (\ref{exp:B}) and (\ref{exp:pseudo-vac}), we obtain
\begin{equation}
    \Phixyz_\lambda(\vec w;v) = 
    \eta^M (\phi_\lambda(\vec w;v) + O(\eta)).
\label{def:phi-lambda}
\end{equation}
It is possible to write down the complicated definition of 
$\phi_\lambda(\vec w;v)$ explicitly, but it is not necessary to our
purpose here. Expanding $\Lambda$ (\ref{def:Lambda}) and
$\Lambda_j^\lambda$ (\ref{def:Lambda-j}) as
\begin{eqnarray}
    \Lambda(u;\vec w; \eta) &=& 
    1 + \eta   \Lambda^{(1)}(u;\vec w) + 
        \eta^2 \Lambda^{(2)}(u;\vec w) + O(\eta^3),
\label{exp:Lambda}
\\
    \Lambda_j^\lambda(u;\vec w; \eta) &=&
    \eta   \Lambda_j^{\lambda (1)}(u;\vec w) + 
    \eta^2 \Lambda_j^{\lambda (2)}(u;\vec w) + O(\eta^3),
\label{exp:Lambda_j}
\end{eqnarray}
we obtain the following equation for the action of leading terms of 
$\hat t(u)$ on $\phi_\lambda(u;\vec w;v)$ from (\ref{t-Phi}) and
(\ref{exp:hat-t}):
\begin{eqnarray}
    \lefteqn{
    \tr\calT^{(2)}(u) \phi_\lambda(\vec w;v) 
    }\nonumber \\
    &=&
    \left( 
      4 \frac{\partial^2}{\partial \lambda^2}
    - 4 \Lambda^{(1)}(u;\vec w) \frac{\partial}{\partial \lambda}
    + 2 \Lambda^{(2)}(u;\vec w)
    \right)
    \phi_\lambda(\vec w;v)
    \nonumber
    \\
    &+& \sum_{j=1}^M 
    \left(
    - 4 \Lambda_j^{\lambda (1)}(u;\vec w) 
        \frac{\partial}{\partial \lambda}
    + 2 \Lambda_j^{\lambda (2)}(u;\vec w)
    \right)
    \phi_\lambda(\vec w_j(u);v).
\label{T(2)-Phi}
\end{eqnarray}

Instead of the Fourier series (\ref{def:Psixyz}), we take a Fourier
transformation of $\phi_\lambda(\vec w;v)$ as a candidate of an
eigenvector of $\tr\calT^{(2)}(u)$ (and hence of $\hat\tau(u)$):
\begin{equation}
    \psi_\nu(\vec w) := 
    \int_{-1}^1 
    e^{\pi i \nu \lambda} \phi_\lambda(\vec w;v) \, 
    d\lambda.
\label{def:Psi}
\end{equation}
Here $\nu$ is an integer to be determined.  Since
$\phi_\lambda(\vec w;v)$ is periodic with respect to
$\lambda$ with period 2, operator $\tr\calT^{(2)}$ acts on
$\psi_\nu$ as follows by virtue of (\ref{T(2)-Phi}):
\begin{eqnarray}
    \tr\calT^{(2)}(u) \psi_\nu(\vec w) &=&
    ( 4 (\pi i \nu)^2 
    + 4 \pi i \nu \Lambda^{(1)}(u;\vec w)
    + 2 \Lambda^{(2)}(u;\vec w))
    \psi_\nu(\vec w) \nonumber \\
    &+& \sum_{j=1}^M ({\rm unwanted})_j,
\label{T(2)-Psi}
\end{eqnarray}
Here ``unwanted terms'' are defined by
\begin{eqnarray}
    \lefteqn{
    ({\rm unwanted})_j = 
    \int_{-1}^1 e^{\pi i \nu \lambda}
    \left(
    - 4 \Lambda_j^{\lambda (1)}(u;\vec w) 
        \frac{\partial}{\partial \lambda}
    + 2 \Lambda_j^{\lambda (2)}(u;\vec w)
    \right)
    \phi_\lambda(\vec w_j(u);v)\, d\lambda}\nonumber
\\
    &=&
    \int_{-1}^1 e^{\pi i \nu \lambda}
    \left( 
    \pi i \nu
    + \sum_{n=1}^N \ell_n
      \frac{\theta_{11}'(w_j - z_n)}{\theta_{11}(w_j - z_n)}
    - \sum_{k=1, k\neq j}^M
      \frac{\theta_{11}'(w_j - w_k)}{\theta_{11}(w_j - w_k)}
    \right) \times
    \nonumber \\
    &\times&
    4 \Lambda_j^{\lambda (1)}(u;\vec w)
    \phi_\lambda(\vec w_j(u);v)\, d\lambda,
\label{unwanted}
\end{eqnarray}
where ${}'$ denotes the derivative. Thus, if $(\dots)$ in
the last expression of (\ref{unwanted}) vanishes, namely, if
\begin{equation}
    \sum_{n=1}^N \ell_n
       \frac{\theta_{11}'(w_j - z_n)}{\theta_{11}(w_j - z_n)}
    =
    - \pi i \nu
    + \sum_{k=1, k\neq j}^M
      \frac{\theta_{11}'(w_j - w_k)}{\theta_{11}(w_j - w_k)}
\label{BA-eq:Gaudin}
\end{equation}
holds for any $j = 1, \dots, M$, then (\ref{T(2)-Psi}) implies that
$\psi_\nu(w_1, \dots, w_M)$ is an eigenvector of $\tr\calT^{(2)}$. 

The corresponding eigenvalue of $\hat\tau(u)$ can be computed from the
expansion (\ref{exp:Lambda}), but a simpler way is the following: Observe
that the eigenvalue of $\tr\calT^{(2)}(u)$ for $\psi_\nu(\vec w)$
$$
      4 (\pi i \nu)^2 
    + 4 \pi i \nu \Lambda^{(1)}(u;\vec w)
    + 2 \Lambda^{(2)}(u;\vec w)
$$
is equal to the coefficient at $\eta^2$ of the right hand side of
(\ref{funct-eq:xyz}) in which we replace $\{ w_j \}$ by a solution to
(\ref{BA-eq:Gaudin}) this time. Thus we can derive the formula for an
eigenvalue of $\hat \tau(u)$ corresponding to $\psi_\nu$ by extracting
terms of order $\eta^2$ in (\ref{q-det-val}) and (\ref{exp:Delta})
(cf.~\cite{skl:87}):
\begin{equation}
    q''(u) - 2 Z(u) q'(u) + (Z(u)^2 - Z(u)') q(u) = \tau(u) q(u),
\label{funct-eq:Gaudin}
\end{equation}
or
\begin{equation}
    \tau(u) = (\chi(u) - Z(u))^2 + \frac{d}{du}(\chi(u) - Z(u)),
\label{tau-val}
\end{equation}
where
$$
    Z(u) :=
    \sum_{n=1}^N \ell_n
       \frac{\theta_{11}'(u - z_n)}{\theta_{11}(u - z_n)}
\qquad
    \chi(u) := \frac{q'(u)}{q(u)} 
    =
    \sum_{m=1}^M \frac{\theta'_{11}(u - w_m)}{\theta_{11}(u - w_m)}.
$$
This result is an elliptic analogue of the spectrum of the XXX Gaudin model
\cite{gau:73}, \cite{gau:76}, \cite{gau:83}, \cite{skl:87}. 

\section{Comments}
We applied the algebraic Bethe Ansatz to the XYZ Gaudin model and
constructed eigenvectors of the generating function of integrals of
motion, $\hat\tau(u)$. Several questions need, however, further
investigation.

First, it would be preferable to have a closed treatment of the
Gaudin model, completely independent of the original XYZ model. The
difficulty is that the expression for $\phi_\lambda(\vec w;v)$ obtained
by differentiation of $\Phixyz_\lambda(\vec w;v)$ is rather voluminous.

Another weak point is that we must assume that the total spin 
$\ell_{\rm total}$
must be an integer. This obstruction comes from the same assumption in
the Bethe Ansatz of the XYZ type spin chains \cite{takh-fad:79},
\cite{take:92}.

It is also an open question if we can obtain all the eigenvectors whithin
our approach. However, this is quite common
situation in the application of the algebraic Bethe Ansatz known as
completeness problem.

We expect that these problems would be overcome within an
alternative approach to the
model known as Separation of Variables, \cite{skl-tak}. The functional
equation (\ref{funct-eq:Gaudin}) should be interpreted
then as a separated equation.

%
%
%
%
\appendix
\renewcommand{\theequation}{\thesection.\arabic{equation}}
\section{Review of the Sklyanin algebra}
\setcounter{equation}{0}
\label{skl-alg}
In this appendix we recall several facts on the Sklyanin algebra and
its representations from \cite{skl:82} and \cite{skl:83}.
We use the notation from \cite{mum} for theta functions:
$$
    \theta_{ab}(z;\tau) = \sum_{n\in\Integer}
       \exp\left( 
             \pi i \left( \frac{a}{2} + n \right)^2 \tau
           +2\pi i \left( \frac{a}{2} + n \right)
                   \left( \frac{b}{2} + z \right)
           \right),
$$
where $\tau$ is a complex number in the upper half plain.

The {\em Sklyanin algebra}, $U_{\tau,\eta}(sl(2))$ is generated
by four generators $S^0$, $S^1$, $S^2$, $S^3$, satisfying 
the following relations:
\begin{equation}
    R_{12}(u-v) L_{01}(u) L_{02}(v) =
    L_{02}(v) L_{01}(u) R_{12}(u-v).
\label{RLL}
\end{equation}
Here $u$, $v$ are complex parameters, the {\em $L$ operator},
$L(u)$, is defined by
\begin{eqnarray}
    L(u) &=& \sum_{a=0}^3 W_a^L(u) S^a \tensor \sigma^a,
\label{def:Lxyz}
\\
    W_0^L(u) 
    &=& \frac{1}{2 \theta_{11}(\eta;\tau)},\qquad
    W_1^L(u)
    = \frac{  \theta_{10}(u;\tau)}
           {2 \theta_{11}(u;\tau) \theta_{10}(\eta;\tau)},
\nonumber
\\
    W_2^L(u)
    &=& \frac{\theta_{00}(u;\tau)}
           {2 \theta_{11}(u;\tau) \theta_{00}(\eta;\tau)},\qquad
    W_3^L(u)
    = \frac{  \theta_{01}(u;\tau)}
           {2 \theta_{11}(u;\tau) \theta_{01}(\eta;\tau)},
\nonumber
\end{eqnarray}
$R(u)$ is {\em Baxter's $R$ matrix} defined by
\begin{equation}
    R(u) = \sum_{a=0}^3 W_a^R(u) \sigma^a \tensor \sigma^a,\qquad
    W_a^R(u) := W_a^L(u + \eta)/W_0^L(u + \eta).
\label{def:R}
\end{equation}
and indices $\{0,1,2\}$ denote the spaces on which operators act
non-trivially: for example,
$$
    R_{12}(u) = 
    \sum_{a=0}^3
    W_a^R(u) 1 \tensor \sigma^a \tensor \sigma^a,\qquad
    L_{01}(u) =
    \sum_{a=0}^3
    W_a^L(u) S^a \tensor \sigma^a \tensor 1.
$$

The above relation (\ref{RLL}) contains $u$ and $v$ as parameters,
but the commutation relations among $S^a$ ($a=0, \dots, 3$) do not
depend on them:
\begin{equation}
    [S^\alpha,S^0    ]_- =
    -i J_{\alpha,\beta} [S^\beta,S^\gamma]_+, \qquad
    [S^\alpha,S^\beta]_- =
                      i [S^0,    S^\gamma]_+,
\label{comm_rel}
\end{equation}
where $(\alpha, \beta, \gamma)$ stands for any cyclic permutation of
(1,2,3), $[A,B]_\pm = AB\pm BA$, and
$J_{\alpha,\beta}=(W_\alpha^2-W_\beta^2)/(W_\gamma^2-W_0^2)$, i.e.,
$$
    J_{12}= \frac{\theta_{01}(\eta;\tau)^2 \theta_{11}(\eta;\tau)^2}
                 {\theta_{00}(\eta;\tau)^2 \theta_{10}(\eta;\tau)^2},
\quad
    J_{23}= \frac{\theta_{10}(\eta;\tau)^2 \theta_{11}(\eta;\tau)^2}
                 {\theta_{00}(\eta;\tau)^2 \theta_{01}(\eta;\tau)^2},
\quad
    J_{31}=-\frac{\theta_{00}(\eta;\tau)^2 \theta_{11}(\eta;\tau)^2}
                 {\theta_{01}(\eta;\tau)^2 \theta_{10}(\eta;\tau)^2}.
$$

The {\em spin $\ell$ representation} of the Sklyanin algebra,
$$
    \rho^{\ell}: U_{\tau,\eta}(sl(2)) 
          \to \End_{\Complex}(\Theta^{4\ell}_{00})
$$
is defined as follows: The representation space is
$$
    \Theta^{4\ell+}_{00} =
    \{f(z) \, |\,
     f(z+1) = f(-z) = f(z), f(z+\tau)=\exp^{-4\ell\pi i(2z+\tau)}f(z) \}.
$$
It is easy to see that dim$\Theta^{4\ell+}_{00} = 2\ell+1$.
The generators of the algebra act on this space as difference operators:
\begin{equation}
    (\rho^\ell(S^a) f)(z) =
    \frac{s_a(z-\ell\eta)f(z+\eta)-s_a(-z-\ell\eta)f(z-\eta)}
         {\theta_{11}(2z;\tau)},
\end{equation}
where
\begin{eqnarray*}
    s_0(z) =  \theta_{11}(\eta;\tau) \theta_{11}(2z;\tau),&\qquad&
    s_1(z) =  \theta_{10}(\eta;\tau) \theta_{10}(2z;\tau),\\
    s_2(z) = i\theta_{00}(\eta;\tau) \theta_{00}(2z;\tau),&\qquad&
    s_3(z) =  \theta_{01}(\eta;\tau) \theta_{01}(2z;\tau).
\end{eqnarray*}
These representations reduce to the usual spin $\ell$ representations of
$U(sl(2))$ for $J_{\alpha\beta} \to 0$ ($\eta \to 0$).
%
%
%
%

\end{document}